\documentclass[a4Paper,11pt]{article}
\usepackage{amsmath}
\usepackage{amsfonts}
\usepackage{amsthm}
\usepackage{enumerate}
\usepackage{enumitem}
\usepackage{latexsym}
\usepackage{pdflscape}
\usepackage{soul}
\usepackage{color,fancybox,graphicx}
\usepackage{lscape}
\usepackage{rotating}
\usepackage{url}
\usepackage{array,multirow}
\usepackage{natbib}
\usepackage[pagewise]{lineno}
\usepackage{tabularx}
\usepackage{caption}
\usepackage{subcaption}

\newcommand{\Dhat}{\hat{D}}
\theoremstyle{definition}

\begin{document}

\newpage
\markboth{The Airport Counter Allocation Problem: A Survey}{A Resource Management
Problem}

\title{\huge \textbf{\Large The Airport Check-in Counter Allocation Problem: A Survey}\\
%~\\
{\normalsize \textbf{T.R.Lalita}\footnote{\emph{Corresponding author}:T R Lalita (trlalita@gmail.com)} \ \textbf{G.S.R.Murthy}\footnote{murthygsr@gmail.com} }\\
{\scriptsize \emph{SQC \& OR Unit, Indian Statistical Institute, St.No.8, Habsiguda,
Hyderabad 500007, India}}}

\maketitle

\abstract{
%Airports are crucial for passenger transport. They are the main mode of
%intercontinental transport. The number of scheduled airline passengers has only
%increased year to year since $2004$, excepting for the closure of air transport
%world-wide due to Covid-19 related lockdowns. As the passenger traffic is set to
%increase again, airports have to brace themselves for lot of passenger related
%congestion at different points in the airport. An important passenger service area
%which has an impact on passenger satisfaction and airport revenue is the check-in
%counters. Decisions made here effect  the movement of passengers in the airport and
%bad decisions can result in chaos. Since $Chun{1996}$, several authors have proposed
%a multitude of models because of variations in the optimization criteria, modeling,
%airport requirements and airport layouts. This paper presents a state-of-art survey
%on the airport check-in counter allocation problem and focuses on relevant models
%and algorithms. We also present a classification of the existing literature based on
%the type of problem solved, objectives considered and methodology considered.
%Particular emphasis is given to research on adjacent scheduling of counters.

An important passenger service area that has an impact on passenger satisfaction
and airport revenue is the check-in counters. The check-in counter allocation
problem consists of allocating adjacent counters to airlines at an airport and 
scheduling the counters through the day subject to operational constraints. It is a
special form of the well-known RCPSP and an NP-hard problem. In addition, the
continuous demand for the highest possible number of counters throughout the day by
each airline makes it a daily challenge for airport operators. As the counters are a
resource, this problem is equivalent to the adjacent resource scheduling problem,
making the solutions for this problem extendable to any adjacent resource
scheduling problem. Decisions made at the check-in counters affect the movement of
passengers in the airport and bad decisions can result in chaos. Since the $1980$s
several authors have proposed a multitude of models with variations in the
optimization criteria, modeling, airport requirements, and airport layouts. This
article presents a state-of-art survey on the airport check-in counter allocation
problem by focusing on relevant models and algorithms. Relevant literature is discussed based on the type of problem solved,
objectives considered and methodology considered. The value of this research is to
help airport operators in planning and allocation of check-in counters, increase airport revenue, improve passenger flow within existing constraints, and optimize utilization of the
existing infrastructure. }

%related to the multiprocessor scheduling problem and the two dimensional strip
%packing problem,
%Classification depends on the solution methodology considered such as simulation,
%linear programming, heuristics and objective.
%We discuss relevant mathematical models, approximation algorithms, heuristic and
%exact approaches from literature
%these papers are also classified according to the main objective minimized in the paper.
%The sources of uncertainty considered are also highlighted.

\section{Introduction} \label{introduction}  % what is the counter allocation problem

%The objective is to assign an optimal number of adjacent counters to airlines. The
%problem is modeled in two parts, the first part determines the minimum counter
%requirement for airlines' daily operation, and the second part is for
%adjacent assignment of counters to each airline. The adjacent allocation problem requires
%all items in the shape of rectangles or irregular structures such as polyominoes, to be
%packed in a minimum area, while the items are not allowed to be moved horizontally and only
%vertical movement is allowed. %ADD THIS TO THE INTRODUCTION.

%Air-travel has become the most preferred mode of travel. The total number of air
%travellers was 4.1 billion in 2017, which is 7.2 per cent higher than the year 2016,
%while the number of flight departures reached 36.7 million in 2017, a 3.1 per cent
%increase compared to 2016 (\cite{ICAO2017}).  The
%present trends in air transport suggest that passenger numbers could double to 8.2
%billion in 2037 (\cite{IATA2018}).

Air-travel has become the most preferred mode of travel. The increase in economic
growth and the average household income has contributed to the growth of the
aviation industry (\cite{zhang2020}). The number of scheduled aircrafts has only
increased year to year since $2004$, excepting for the closure of air transport
world-wide due to Covid-19 related lockdowns. As the passenger traffic is set to
increase again, airports have to brace themselves for a lot of passenger related
congestion at different points in the airport. Check-in counter allocation is an
important resource scheduling problem at airports. The ever increasing demand for
airport resources has made careful planning and optimal use of resources a
necessity. Delay due to inadequacy or inefficient management of airport facilities
may result in penalty for the airline companies (\cite{hsu2005scheduling}). It has
been found that 80\% of the passenger delay at airports is due to waiting for
check-in (\cite{takakuwa2003modeling}, \cite{parlar2013allocation}). To overcome
this problem, kiosks were introduced in airports. \cite{abdelaziz2010study} study
outcomes of introducing kiosks at Cairo airport. Passengers with baggage were asked
to use luggage-drop counters after check-in at kiosks. This reduced waiting time at
kiosks for check-in, but resulted in waiting at luggage drop areas. Though airlines
rely on kiosks for management of queue and congestion at airports, their use mainly
depends on the convenience for check-in as well as luggage drop. Adding additional counters at the airport is not generally a choice as it requires area in the airport, conveyor belts and also interlinking of the existing multilevel system of baggage transit with the new belt.

%Even though many
%passengers check-in to their flights online, luggage drop areas are still expected
%to be congested()
% add this following paragraph somewhere

Therefore, efficient usage of check-in counters will improve service levels at airports,
result in smaller queue lengths for airport operators and faster check-ins for
passengers. As most of the passenger waiting time is spent in the check-in area at an
airport, a reduction impacts public perception of the level of service and indirectly
enhances the airport revenue due to the increased stress free time in the commercial
areas of the airport (\cite{hsu2005scheduling}, \cite{lin2013passengers},
\cite{parlar2013allocation}).

To achieve all the above objectives in counter allocation, different optimization
techniques have been used in the literature. This paper reviews all the methods proposed
for counter allocation. The techniques discussed in this paper for allocating adjacent
counters can be applied to other adjacent resource allocation problems, such as:
warehouse space optimization where warehouse area has to be assigned to customers for
storage for a certain time, berthing problem at ship yards, where berthing space has to
be allocated to ships given the ship size and time for which it docks at the shipyard
(see \cite{bierwirth2010survey}), for assignment of computer hard disk memory (if
contiguous allocation of memory is required (see \cite{Duin2006})) and other resource
scheduling problems where jobs cannot be shifted in time but resource requirements may be
satisfied by any set of adjacent resources and may vary with time.

Efforts at reviewing resource usage at airport terminals have not explicitly focused on
check-in counter utilization. \cite{cheng2012theory} and \cite{wu2013review} have
reviewed models presented in literature for resource scheduling at airport terminals.
\cite{cheng2012theory} reviews theory of allocating and scheduling resources by grouping
existing literature into three parts, viz., methods of integrating airport operation
data, methods of predicting passenger flow at airport terminals, and optimization of
allocating and scheduling passenger service resources at airport terminals.
\cite{wu2013review} review airport passenger terminal models by classifying them based on
usage scenarios. \cite{wu2013review} review work related to capacity planning,
operational planning and design, security policy and airport performance.
\cite{cheng2012theory} and \cite{wu2013review} do not study the check-in counter
allocation problem in airport terminals. In this article, the authors provide technical guidance to the reader on how to evaluate the existing methods on counter allocation and appropriately model. This article simplifies the usage of existing knowledge by proposing a unified framework for modelling the counter allocation problem. Under this framework, we scrutinize different works
based on their applicability, constraints, objectives considered and the time taken
to solve the problems. This work gives any airport operator the ability to navigate
through the articles published since four decades ago and implement a suitable model.
% Our paper enables the reader to pursue modelling his check-in counter allocation
% problem,  based on his requirements.

This article has five sections. Section 2 describes the check-in counter allocation
problem. Section 3 classifies different methods in the literature by the problem type
addressed and its limitations. Section 4 discusses different approaches for modelling the
check-in counter allocation problem. Section 5 gives the conclusion.
%This paper reviews published literature related
%to planning and modelling the counter allocation problem for flight departures at an
%airport, it does not consider studies concerning with real-time management of
%sudden/unforeseen crisis/circumstances at the airport.

\section{The Check-In Counter Allocation Problem} \label{theproblem}

\begin{figure}
  \centering
  \includegraphics[scale=0.6]{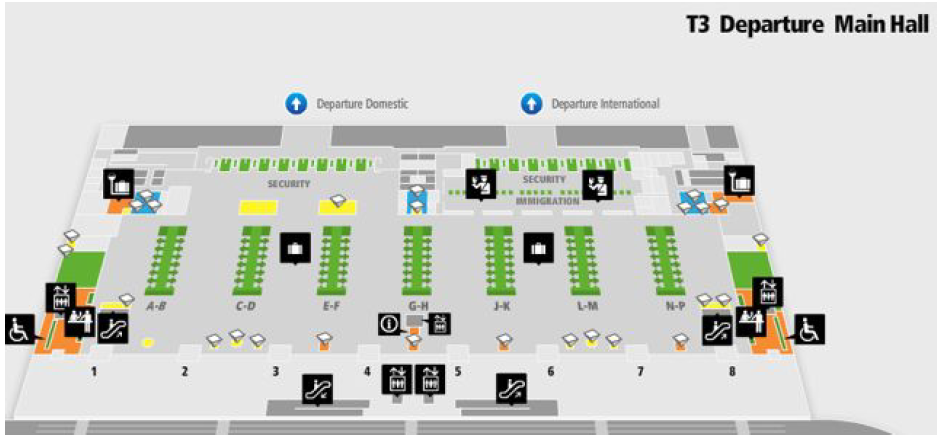}
  \caption{Departure Hall in an Airport}\label{DepHall}
\end{figure}

Airports are divided into airside and landside areas. Airside operations consist of
scheduling flights on the runways and other related on-air flight operations. Landside
operations consist of scheduling and organization of the processes that passengers need
to undergo before boarding a flight. It is imperative that airport landside operations
focus on timely boarding of passengers and flight departures. For ease of operation, an
airport has designated areas for different pre-boarding operations for passengers, such
as check-in, security check, immigration, etc. In most airports, the check-in area
consists of multiple structures of counters arranged around a conveyor belt in a u-shape.
This structure is referred to as an island. Many such islands make up the check-in area
(see Fig:~\ref{DepHall}).

\begin{figure}
  \centering
  \includegraphics[scale=0.4]{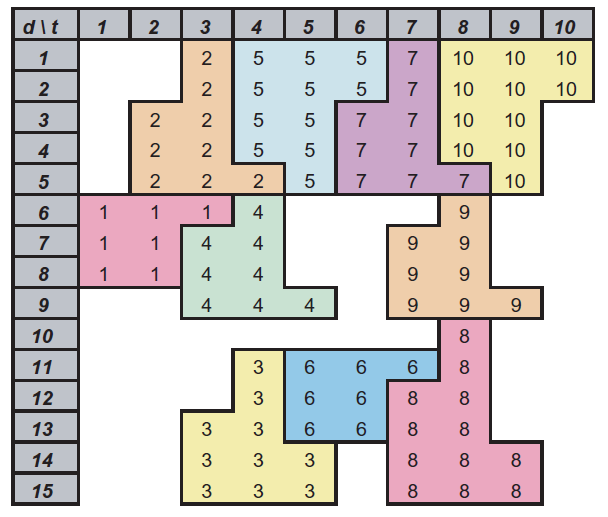}
  \caption{Polyominoes in counter allocation }\label{VanDijk1}
\end{figure}

Though all the counters are not physically adjacent to one another, this is supposed to
be the case for simplicity in modelling the problem. %and computation.
Moreover, counters can be rearranged according to physical layout at later stages for
implementation of the solution. In most of the literature discussed in this paper, a two
dimensional space is used as a representation of the counters in the planning horizon
(see Fig.\ref{VanDijk1}). The planning horizon is the time period for which counter
assignments are computed. To model the problem, the planning horizon is divided into
smaller Time Intervals (also Time Windows(TWs)). The length of the TW has to be chosen
while modelling a problem. Too small TWs mean less time for staff to handle operations
and large TWs can consist of large variations in passenger arrivals. Fig:~\ref{VanDijk1}
shows counter allocation by \cite{dijk2006} with a TW length of 60 min, planning horizon
of 10 hours with resource availability of 15 counters. Here, the numbers represent the
flights, `d' represents the desk number or counter number and `t' represents TWs. Since
resource requirement times are fixed, it is not possible to shift items horizontally,
only vertical movement is allowed.

\begin{figure}
  \centering
  \includegraphics[scale=0.5]{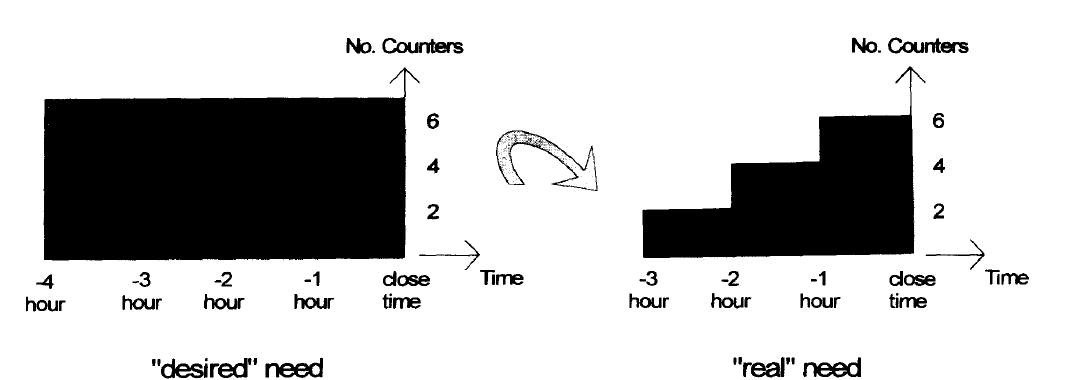}
  \caption{Desired need for airlines as compared to real need}\label{RealNeed}
\end{figure}

\begin{figure}
  \centering
  \includegraphics[scale=0.5]{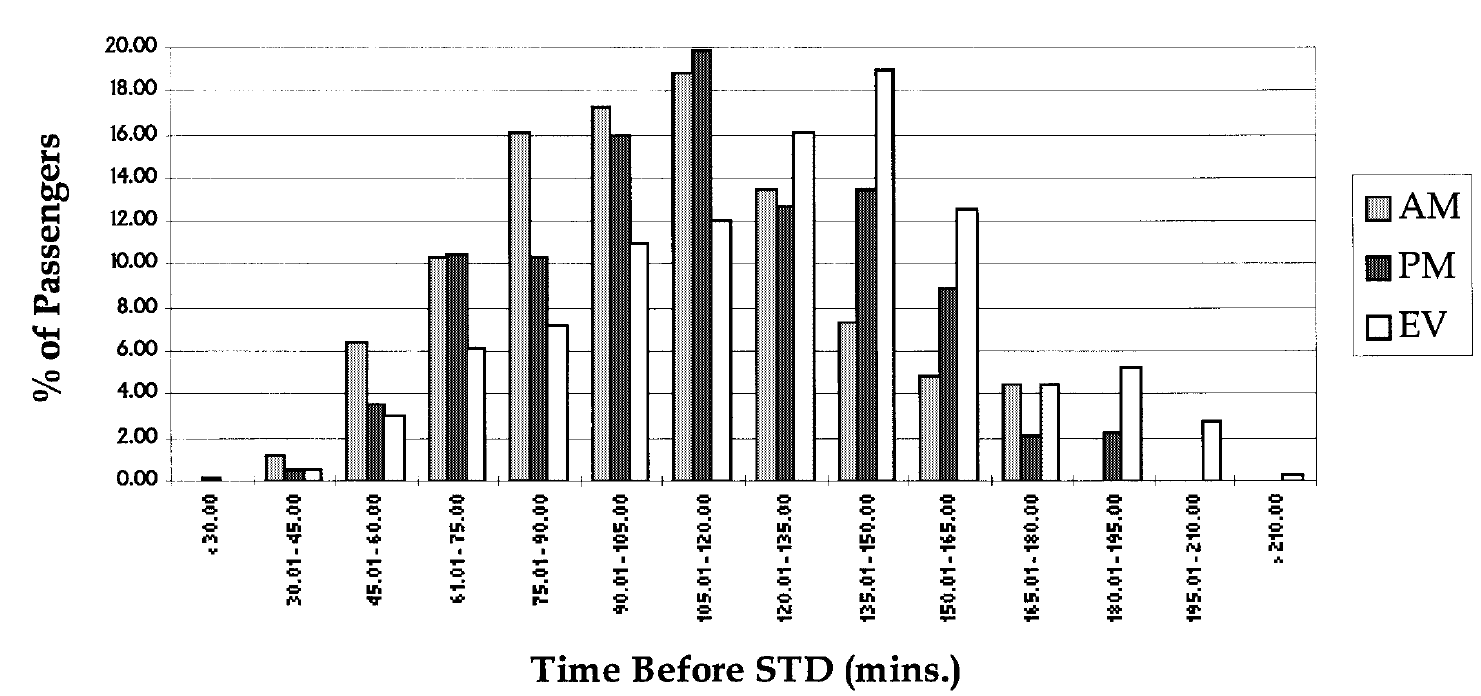}
  \caption{Arrival Distribution Observed at }\label{ChunArrivalDist}
\end{figure}

\begin{figure}
  \centering
  \includegraphics[scale =0.5]{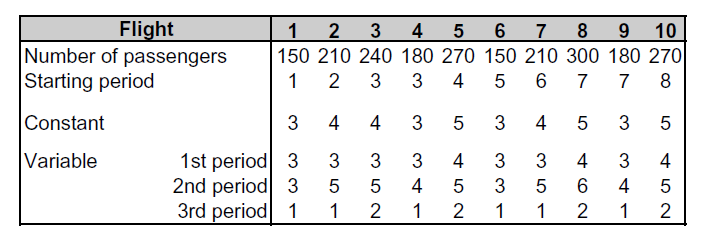}
  \caption{Minimum Counters Required in constant and variable case}\label{CounterReq}
\end{figure}

These counters are owned by the airport operator and suitably leased to airlines. The
airport operator has to issue the minimum number of adjacent counters required for each
flight/ group of flights. Adjacent counters are important for passenger and airline
convenience, airline visibility and operations, and for ease in baggage management and
sorting, since baggage from adjacent counters of an island is collected at one baggage
collection center. Typically, airlines demand more counters from the airport operator
than required, for ease in providing service and visibility at the airport (see Fig.
~\ref{RealNeed} from \cite{chun1996scheduling}). This creates problems for the airport
operator, who has to now simultaneously satisfy the real need for counters, ensure
optimal allocation to all airlines and minimize changeover operations. A changeover
operation consists of staff of one airline closing the counters at the end of a TW and
shifting to other counters for continuing the check-in process in subsequent TW(s). To
minimize the changeover operations, counters were traditionally allocated to a flight for
its entire check-in time. Since this type of allocation would result in a waste of
counter time due to time-varying arrival distribution, time-varying counter requirements
are proposed (\cite{dijk2006}, \cite{bruno2010mathematical}, \cite{araujo2015optimizing}
, \cite{opt2019} etc). This results in structures known as polyominoes (see Fig.:
\ref{VanDijk1} for each flight (or a group of flights). Actual counter requirement (for
illustration see Fig. ~\ref{CounterReq}, from \cite{dijk2006}) is computed using the
arrival distribution of passengers (for illustration see Fig.:~\ref{ChunArrivalDist} from
\cite{chun1999intelligent}, where arrival distribution is computed for flights in the
morning, afternoon and evening), and an appropriate load factor. Load factor is defined
as the percentage of passengers of an airplane seat capacity expected to finally board
the flight. The counter allocation problem is reduced to determining counter requirement
of flights and then placing the resulting polyominoes, that can be moved only vertically,
in the counter-TW area in an optimal way.

\begin{figure}
\centering
\begin{subfigure}{.49\textwidth}
  \centering
  \includegraphics[scale=.7]{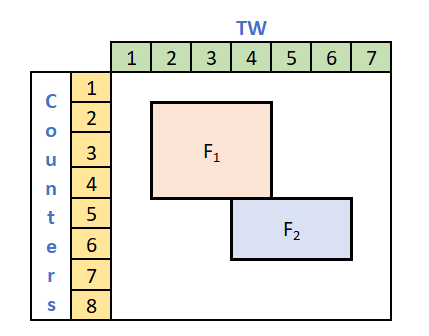}
  \caption{Static Counter Allocation}
  \label{fig:static1}
\end{subfigure}
\begin{subfigure}{.49\textwidth}
  \centering
  \includegraphics[scale=.7]{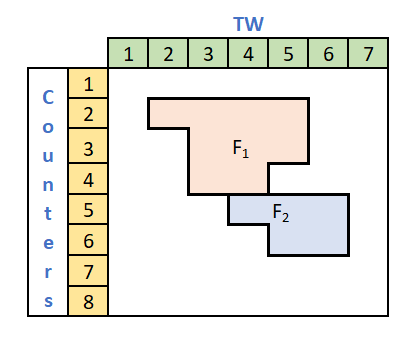}
  \caption{Dynamic Counter Allocation}
  \label{fig:static2}
\end{subfigure}
\caption{Dynamic Counter Allocation saves counter time}
\label{StaticDynamic}
\end{figure}

It has been observed that time-varying or dynamic counter allocation
(Fig.\ref{StaticDynamic}) provides on-time service and consequently results in
significant savings in terms of counter time, operation costs and reduced queue size
(\cite{joustra2001simulation}, \cite{dijk2006}). In contrast, static counter allocation
assumes constant counters in the check-in period. In both types of allocation, two or
more flights of an airline can be grouped together for assignment if the demand for
counters overlaps in some TWs, i.e., if flights are scheduled for departure a few hours
apart. It is general practice by airlines to assign common counters to a set of flights
with simultaneous demand in at least one TW. The demand for counters is treated as demand
for one departure. This enables passengers of an airline to use the same counters
irrespective of the flight boarded. Dedicated counters on the other hand exclusively
serve passengers of a single flight. Both allocations are used in airports based on
opportunities available to group flights. This allows the airlines to minimize the
counter operating cost and changeover operations also. In most of the models proposed for
counter determination, minimizing the counter operating cost is the objective. Counter
operating cost includes the cost of operating baggage belts, the cost of staff operating
the check-in counters, queue cost and the cost of delay. \cite{hsu2005scheduling}
construct various cost functions related to facility management.

 Some of the challenges to this problem, as discussed by \cite{snowdon1998avoiding},
are estimating the passenger arrival distribution, the complexity of determining the
passenger mix (number of passengers using different services), changes to the resources
(mainly check-in facilities) available depending on the arrival distribution, ensuring
passenger service levels are attained and evaluating flight data in order to group
flights which can benefit from common counters vis-a-vis dedicated counters. Initial
attempts to solve the problem through simulation improved the prevailing methods for
counter allocation. \cite{atkins2003right} propose simulation to compare operational
strategies and to determine the optimal staff levels required. \cite{chun1996scheduling}
and \cite{chun1999intelligent} use simulation for counter determination and counter
allocation. Real life modelling of the check-in process at airports has been presented by
\cite{joustra2001simulation} and \cite{dijk2006} for check-in at Schiphol airport in
Amsterdam, \cite{atkins2003right} at Vancouver airport, \cite{bruno2010mathematical} at
Naples International Airport, \cite{lous2011modelling} at the Copenhagen airport,
\cite{al2016optimization} at an airport in Kuwait and \cite{felix2017hybrid} at the
airport of Lisbon and \cite{chun1996scheduling}, \cite{chun1999intelligent} at the Hong
Kong airport etc. The counter allocation problem with the adjacency restriction is
NP-complete (\cite{dijk2006}, \cite{Duin2006}) and cannot be solved in polynomial time.
The complexity of the problem has been studied in detail by \cite{Duin2006}. All the
models proposed for counter allocation including the above real-world applications have
been classified on the basis of the problem solved. These methods are discussed in detail
below.

\subsection{Determining Optimal Number of Check-in Counters}  \label{CounterDetermination}

This section discusses the problem of determining counter requirement for a flight or a
set of consecutive flights of an airline. The number of counters allocated to a flight
(or group of flights) depends on the arrival pattern of passengers, the queueing area
available, queue length in an interval and restrictions on waiting time. Different
authors have modelled the problem with different constraints, different objectives and
different facilities (eg: counters and kiosks). Before studying these models, we present
a basic model for counter determination:

\begin{align}
\texttt{Min~~} & \sum_{ij} x_{ij} \label{basic_obj} \\
\texttt{subject to} & \\
& \sum_{j} a_{ij}s_{i} \leq \sum_{j} x_{ij}t  \label{basic1} \\
& x_{ij} \geq 0 \texttt{ are nonnegative integers.}
\end{align}

\begin{figure}
\centering
\begin{subfigure}{.49\textwidth}
  \centering
  \includegraphics[scale=.6]{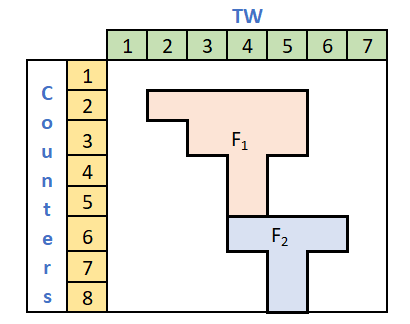}
  \caption{Allocation using Basic Formulation}
  \label{fig:sub1}
\end{subfigure}
\begin{subfigure}{.49\textwidth}
  \centering
  \includegraphics[scale=.6]{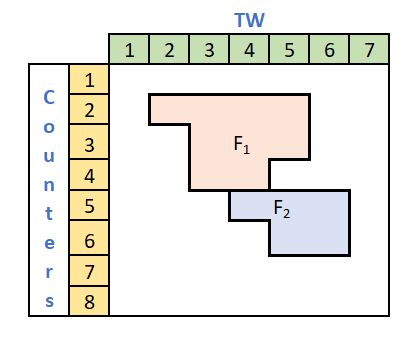}
  \caption{Improved Allocation}
  \label{fig:sub2}
\end{subfigure}
%\label{improvement}
\caption{Reduced Changeover Operations}
\end{figure}

In this formulation, $x_{ij}$ is the number of counters assigned to the $i^{th}$ flight
in the $j^{th}$ TW, $a_{ij}$ is the average number of arrivals for $i^{th}$ flight in
$j^{th}$ TW (this is obtained from passenger surveys conducted at the airport), $s_i$ is
the average service time per passenger for the $i^{th}$ flight, `t' is the length of each
TW. This formulation provides resources exactly in proportion to airline requirement.
This results in counter allocation with peaks and troughs exactly like the passenger
arrival distribution. Since airport operators mandate airlines to limit waiting times and
counter queue lengths, the only way to improve the solution to this model and get more
rectangular polyominoes is to postpone the check-in of some passengers while respecting
the service level requirements. Note that with even a slightly more rectangular structure
of the counters allocated (see Fig.(\ref{fig:sub2})), changeover operations are reduced,
improving the basic solution (see Fig.(\ref{fig:sub1})).

\subsubsection{Mathematical Models for Counter Determination}\label{SectionCDetermination}
Published work on check-in counter determination is described briefly in the following
paragraphs.

\cite{park2003optimal} published a paper for optimal assignment of check-in counters.
Their paper aims to assign counters based on passenger arrival distribution at the
airport. Other factors considered are aircraft type (standard or chartered), aircraft
size, time allowed for check-in, passenger arrival distribution, ticket status (economy,
business, first class etc), processing time of staff at the check-in counters and load
factor assumed. A passenger survey at the airport (in \cite{park2003optimal}) determines
the variation in load factors during peak and non-peak hours and the resulting arrival
patterns. The airport arrival patterns are then used as input to a regression model to
determine the cumulative arrival pattern based on time before departure. Counters are
allocated to airlines directly in proportion to the estimated passenger arrivals. This
kind of allocation may not result in an optimal assignment of counters to airlines since
the overall cost to the airport operator or queue lengths among other things are not
considered.

\cite{bruno2010mathematical} propose the following static model to determine the optimal
number of counters for flights with the objective of reducing counter cost and queue
length.

\begin{align}
 \texttt{ Minimize~~ }z & = \sum_{j} \sum_{t} (h_j.I_{jt} +s_j.x_{jt})   \label{obj_bruno} \\
\texttt{  subject to } & \nonumber \\
    I_{jt} & = I_{j(t-1)} + d_{jt} - q_{jt}, \ j = 1,2,...,J, \ t = 1,2,...,N, \label{bruno1} \\
    p_jq_{jt} & = C_tx_{jt}, \ j = 1,2,...,J, \ t = 1,2,...,N, \label{bruno2} \\
    \sum_{j} p_jq_{jt} & \leq C_t, \ t = 1,2,...,N, \label{bruno3} \\
    I_{jt} & = 0, \ t \in \ T_j \label{bruno4} \\
    q_{jt}, \ I_{jt} & \geq 0, \ j = 1,2,...,J, \ t = 1,2,...,N, \label{bruno5} \\
    x_{jt} & \in {0,1}, j = 1,2,...,J, t = 1,2,...,N, \label{bruno6}
\end{align}

The following notations are used in the above model. $h_j$ is the cost associated with
queue related to flight $j$, $s_j$ is the desk opening cost for flight $j$, $T$ is the
planning horizon (usually one day), $l$ is the length of the TWs considered, $N$ is the
number of TWs, $J$ is the number of flights scheduled in $T$, $p_j$ is the average desk
service time for flight $j$, $d_{jt}$ is the service demand for flight $j$ in TW $t$,
$C_t$ is the available check-in time based on counters operating in TW $t$, $I_{j0}$ is
the number of passengers of flight $j$ waiting before counters open for flight $j$, $T_j$
is the set of TWs in which counters for flight $j$ do not operate. Decision variables
are: $I_{jt}$, the number of passengers in queue for flight $j$ at the end of TW $t$,
$q_{jt}$, the number of passengers of flight $j$ to be accepted for service in TW $t$.
Even though $x_{jt}$, the binary variable representing the possibility of checking-in
passengers for flight $j$ in TW $t$, is defined as a decision variable, it is not, since
this is fixed in advance and cannot be restructured. Constraints (\ref{bruno1}) represent
the change in queue length between two successive TWs. Constraints (\ref{bruno2}) ensure
that enough counter time is available for passenger check-in in TWs where check-in is
possible. Constraints (\ref{bruno3}) ensure the overall service capacity is as required,
constraints (\ref{bruno4}) forces all passengers of flight $j$ to be accepted by the
closing time of check-in service for flight $j$.

\cite{bruno2010mathematical} propose models for both static (see model
(\ref{obj_bruno})-(\ref{bruno6})) and dynamic counter allocation. Both the models define
the number of passengers in queue for a flight and the passengers accepted for service in
each TW as decision variables. In the static model presented above, the total cost of
counter operation and the cost of queue is minimised (see objective function
(\ref{obj_bruno})). In the dynamic model, counters operating in each TW and the cost
associated with queue are minimized. The authors derive mathematical formulations from
the Capacitated Lot Sizing problem in literature (see \cite{bitran1982computational}) and
\cite{florian1980deterministic}). The authors also present a real-life airport management
study at the Naples airport.

The model presented by \cite{araujo2015optimizing} is an extension of the model by
\cite{bruno2010mathematical}. \cite{araujo2015optimizing} present two models for
determining counter requirement for flights and aim to determine the optimal number of
counters for flights operating at an airport. ILPs are presented for dedicated and common
counter check-in. Following is the ILP for dedicated counter allocation.

\begin{align}
 \texttt{ Minimize~~ }z & = \sum_{j} \sum_{t} (h_j.I_{jt} +s_j.x_{jt})   \label{ar_obj} \\
\texttt{  subject to } & \nonumber \\
    I_{jt} & = I_{j(t-1)} + d_{jt} - q_{jt}, \ j = 1,2,...,J, \ t = 1,2,...,N, \label{ar1} \\
    \sum_{j} p_jq_{jt} & \leq C_t, \ t = 1,2,...,N, \label{ar3} \\
    I_{jt} & = 0, \ t \in \ T_j \label{ar4} \\
    I_{jt} & \leq \alpha.(d_{jt} + I0_{jt}), \ j = 1,2,...,J, \ t = 1,2,...,N, \label{aandr6} \\
    q_{jt}, \ I_{jt} & \geq 0, \ j = 1,2,...,J, \ t = 1,2,...,N, \label{ar5}
\end{align}

An additional service level constraint (such as constraint (\ref{aandr6}) for dedicated
counter allocation as above) is added to both models to ensure that only a small
percentage of passengers, $\alpha$, remain in queue at the end of a TW. The service level
is chosen by the airport operator. The models predict the counter requirement for flights
(or for a group of flights). The model solutions are analysed using simulation for model
validation in different scenarios that may occur at the airport. Once the models are
validated, the corresponding solution is implemented by allocating adjacent counters to
each flight using the ILP by \cite{dijk2006}. Though a service level constraint is added
to the model, the waiting time of passengers and the maximum queue length allowed are not
considered. It can be easily proved that model solution does not follow FIFO and due to
this counter usage may be different than expected.

\cite{opt2019} propose a model ((\ref{objUVF})-(\ref{con4UVF})) for determining the
optimal number of counters with postponement of service times. The model considers the
passenger arrival distribution, where $d_{ji}$ is the arrivals in TW $i$ for flight $j$,
$u_{ji}$ is the number of arrivals (in TW $i$), served in TW $i$ and $v_{ji}$ is the
number of passengers arrivals for flight $j$ served in TW $i+1$. $P$ denotes the time
horizon for these departures. Let $i_o(1)$ denote the counter opening time of the first
departure and $i_c(1)$ denote its closing time, then, $P
=\{i_o(1),i_o(1)+1,\ldots,i_c(\Dhat)\}$, $P_j={i_o(j),i_o(j)+1,\ldots,i_c(j)}$ and
$d_{ji}=0$ for all $(j,i)$ such that $i\not\in P_j$.

Constraints (\ref{con1UVF}) and (\ref{con2UVF}) ensure that all passengers are served,
constraint (\ref{con3UVF}) ensures that enough counter time is available for serving all
the passengers of all the flights. The formulation makes the resulting counter allocation
more rectangular by postponing service to the next TW for some passengers. It also limits
the waiting time by length of two TWs and follows the FIFO queue discipline. Adding a
service level constraint to this model may improve counter allocation. \cite{opt2019}
also propose a model for adjacent check-in counter assignment, discussed in section
\ref{CounterAllocation}.

\begin{align}
    \texttt{ Minimize~~ }& z \ \label{objUVF}  \\
     \texttt{  subject to } & \nonumber \\
     u_{ji} + v_{ji} & = d_{ji}, \forall \ i \in P_j, \ j=1,2, \ldots, \Dhat, \label{con1UVF}  \\
     v_{ji_c(j)} & = 0\  \forall \  j=1,2,\ldots, \Dhat,  \label{con2UVF} \\
     \sum_js_j(u_{ji} + v_{j(i-1)}) & \leq \omega c_i\ \forall\ i \in P, \label{con3UVF} \\
     c_i & \leq z\ \forall\ i \in P, \label{con4UVF} \\
     u_{ji},  \ v_{ji} \texttt{ and }  & c_i \ \forall\ i \in P,\texttt{ are nonnegative integers.}
\end{align}

\cite{hsu2012dynamic} also propose an ILP to determine the number of check-in facilities
required for a departure from the airport. Check-in facilities considered are the
counters and kiosks, online check-in, and barcode check-in. Check-in services offered by
these facilities are : ticket purchase, check-in, boarding pass and checking baggage,
each offering one or more services. Each facility may offer different check-in services.
The model proposed aims to minimize the passenger waiting time, operation costs and
counter requirements. The authors explore dynamic allocation of different check-in
facilities and passengers at the airport. For modelling the problem only two types of
check-in facilities, counters and kiosks are considered. Passengers are assigned
different facilities based on service requirement. The model proposed is based on the
dynamic model by \cite{nikolaev2007sequential}. The services required by passengers and
their arrival times are predicted and are an important input to the model. The assignment
of the $n^{th}$ passenger by the model depends on the assignment of the $(n-1)^{th}$
passenger. Therefore, for dynamic assignment of passengers, various possible scenarios
are analysed. Due to this, the model takes more than 3 hours to solve for 15 passenger
arrivals. In view of this the authors have used clustering algorithms (in heuristics) for
solving the model. A real life scenario has been presented in the study to show improved
counter utilization rates on using the model at an airport in Taiwan. A brief analysis of
the modelling of the problem suggests that the success of the model is entirely dependent
on the willingness of the passengers to use the check-in facility assigned to them. The
model may fail to generate an effective solution at airports where passengers prefer to
use other facilities, it may cause the check-in facilities operating at the airport at
certain TWs to be insufficient.

Determining the counter requirement at an airport is necessary but sometimes, the
capacity of check-in counters at an airport also needs to be assessed. It is essential in
decision making for airport expansion. Airport terminal capacity is assessed by
\cite{brunetta1999operations}, \cite{fayez2008managing} and \cite{pacheco2003managerial}.
\cite{kiyildi2008capacity} present a fuzzy logic method to determine check-in capacity at
Antalya airport in Turkey.

Capacity of a check-in counter is computed as the number of passengers and luggage that
can be checked-in in one hour. In the model, possible passenger and luggage combinations
are determined. Number of passengers travelling together and the number of bags (luggage
carried by a person/group) is randomly generated and fuzzy logic is used to predict the
processing time. It was observed that the processing time is highly correlated with
volume of luggage (high $R^2$). Fuzzy models are used to calculate the capacity of each
check-in counter in terms of the number of persons served and the number of luggage items
checked-in. This is used to calculate the check-in capacity of the airport which is
compared with the current passenger arrivals and growth rate at the airport for adding
new resources to the airport.

\cite{hwang2012airport} present a mathematical model for optimizing check-in counters,
kiosks, part-time staff and full-time staff required during each shift in a week at
airports. Their model assumes static counter allocation and defines the number of
passengers using counters and kiosks in a shift as parameters. The study conducted
gathers information on counter requirements and cost of operation in the presence of
various factors including the day of the week and varying load factors. The authors
compute the ratio of counters to kiosks that would be best suited for serving a flight.
Model solutions are verified for optimizing cost in different scenarios using simulation.
The authors conclude that usage of kiosks reduces the operational costs and that services
can be improved by installing additional kiosks. Regarding the average waiting time, the
study concludes that check-in counter cost can be reduced by limiting passenger service
time. Though the observation is fact-based, it is not clear as to how this may be
achieved. This recommendation also contradicts the existing counter planning strategies
of using average service time to determine counter allocation instead of fixing service
time per passenger. Reducing service time per passenger would definitely make the counter
allocation problem easier, but it may not be always possible to limit passenger
processing time.

Simulation has been used for counter determination by \cite{chun1996scheduling},
\cite{chun1999intelligent} and \cite{dijk2006}.

\subsection{Adjacent Counter Allocation} \label{CounterAllocation}
The most challenging problem at an airport is adjacent counter allocation for a flight
while simultaneously maximizing counter utilization. \cite{chun1996scheduling} first
addressed this problem by defining structures called counter profiles.

\begin{figure}
  \centering
  \includegraphics[scale=0.5]{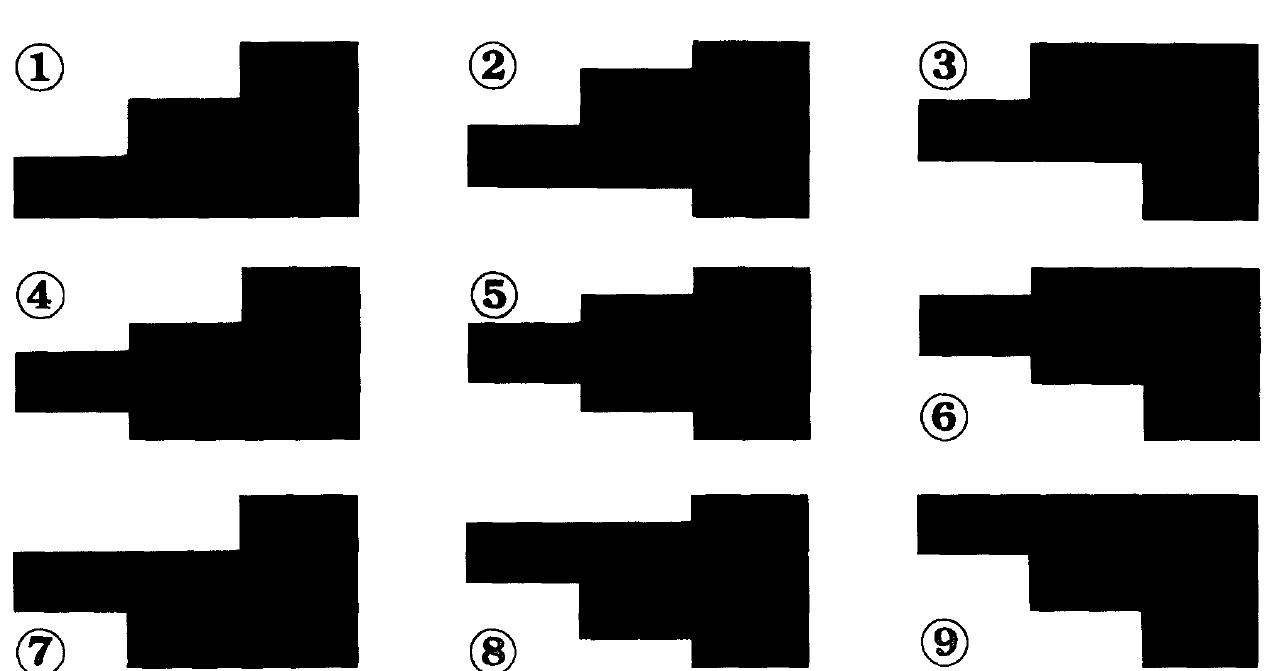}
  \caption{Possible variations of the 2-4-6 Counter Profile} \label{ChunCounterGlob}
\end{figure}

\begin{figure}
  \centering
  \includegraphics[scale =0.5]{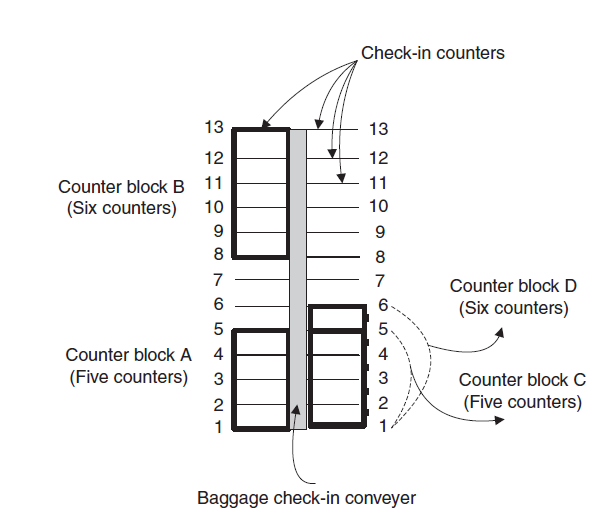}
  \caption{Blocking by \cite{yan2004model}}\label{CounterBlocks}
\end{figure}

A counter profile is a two dimensional shape that defines the check-in counter resource
requirement for a flight. In the algorithm proposed by \cite{chun1996scheduling},
different operators are presented to change the shape of the counter profile to check for
a convenient fit for allocation in the counter-TW rectangle of fixed dimensions. All
possible shapes are considered for every counter profile (see Fig.\ref{ChunCounterGlob}).
Placing a counter profile in a two dimensional counter TW space ensures adjacency of
counters. Algorithm by \cite{chun1996scheduling} starts with selecting one possible
counter profile for a flight, it is assigned a location in the counter TW space (also
Gantt Chart), all the relevant constraints are checked, if some of the constraints are
violated (i.e. if the specific shape of the counter profile cannot fit in the Gantt
Chart), the algorithm backtracks by removing the previous allocation to a flight and
freeing up space. The counter profile shapes may need to be changed in order to
accommodate the remaining flights. \cite{chun1996scheduling} simulated different possible
shapes of the counter profile, similar to fitting pieces of a puzzle randomly until an
acceptable solution was reached. If the order in which flight counter profiles are
allocated does not have a feasible solution, the algorithm involves backtracking to
remove flight allocations. This results in a time consuming hit and trial process for
counter allocation. An ILP is presented by \cite{dijk2006} for the same problem. The ILP
considers all possible arrangements of a counter profile in the counter-TW area, till an
optimal solution is obtained. Due to this, the ILP takes a long time to converge to a
feasible solution (especially for large flight departures). The formulation is effective
for small size problems compared to the real-world problems of the day. In the static
counter allocation problem considered by \cite{yan2004model} and
\cite{yan2005minimizing}, each counter has one or two service lines to provide check-in
service to passengers. To allocate adjacent counters, the authors define each block as a
set of service lines. \cite{tang2010network} defines each block as a set of adjacent
counters. Two blocks may overlap as a service line/ counter may belong to both (see
Fig.~\ref{CounterBlocks}).

Service lines are demanded by airlines in accordance with the passenger arrival pattern
and the number of passengers on the flight. \cite{yan2004model} propose a static model to
assign blocks to flights. The objectives of the study are to allocate counters to
minimize passenger walking distance and reduce inconsistency in counter location.
Allocation to different blocks of counters on different days of the week for the same
flight is defined as an inconsistency. Inconsistency values are introduced as the
airlines prefer to have the same block of counters allocated for flights in a week.

Passenger walking distance is calculated for possible flight assignments to a block and
all possible blocks to which a flight can be assigned are considered in the model. For a
given flight and block combination, it is the average distance a passenger walks after
check-in from the block allocated till boarding. An ILP is proposed for minimizing the
passenger walking distance subject to allowable inconsistency. The constraints of the ILP
need preprocessing to exclude redundant constraints. Since the solution method lacks
scalability, three heuristic models are proposed to solve the problem. The first model
provides the assignment for a day. Based on this, inconsistency values are set and the
second model is solved for a minimum inconsistency value. Then, the third model is solved
for the final assignment of flights to blocks. A heuristic is thus provided for solving
the model for a single day. This heuristic has been used in real-life modelling at an
airport in Taiwan. The models are used for allocating 140 counters to 70 flights
departing from the airport. The number of variables exceed 80000 and constraints exceed
70000. Results from the case study conducted for Taiwan Airport show reduction in the
passenger walking distance by 4\%.

\cite{yan2005minimizing} study the dynamic counter allocation problem. Their objective is
to allocate adjacent counters to flights by minimizing the total inconsistency in blocks
allocated to a flight. Blocking plans are similar to blocking by \cite{yan2004model}.
Each counter may have many service lines (as also defined by \cite{yan2004model}) and the
service line requirement for each flight in each TW is known. Since the number of service
lines differ from one block to another, additional lines need to be opened or closed
between successive TWs. This adjustment of service lines between two consecutive TWs is
considered an inconsistency. From the inconsistency values attached to blocks pairwise,
the inconsistency value in flight allocation on a day is computed. The problem was solved
by an ILP, allocating a block (possibly different) to each flight in each TW. Due to the
complexity involved in computing the inconsistency values, the problem becomes
increasingly complex with increase in number of flights, and the authors propose a
heuristic algorithm. The heuristic for the problem divides it into subproblems, each of
which is solved. Each subproblem comprises of equal counters and flight departures. An
exchange mechanism between the two parts is proposed to proceed further (similar to the
evolutionary algorithm by \cite{mota2015check}). The exchange of two flights, one
selected from each group, (with overlapping departure times), and then re-solving for a
solution is continued till a reduction in inconsistency is achieved. The allowable
inconsistency value is chosen as appropriate. The solution with the least inconsistency
value (according to the stopping criterion) is chosen. The order of exchanges effects the
final solution, hence, further improvements to the heuristic algorithm are suggested. For
large problems, dividing flights into two groups may be insufficient. Due to this number
of groups is increased and the solution algorithm is improved. A disadvantage is the
increase in complexity of the problem with additional subgroups, as a result of increase
in the number of flights.

\subsubsection{Mathematical Models to ensure Adjacency} \label{SectionCAllocation}
\cite{dijk2006} model the problem combining both simulation and integer programming. The
objectives of the study are to determine the minimal counter requirement for each flight
and then to allocate counters at the airport with adjacency maintained. The problem is
solved in two stages. In the first stage terminating simulations are run to determine
counter requirement till the solutions satisfy service level requirements. In the second
stage, an ILP is solved for adjacent allocation of counters. The study presents ILPs for
both variable and constant counter allocation. These models ensure adjacency of the
counters allocated to each flight and solve to an optimal solution. The ILP for variable
counter allocation is given below:

\begin{align}
 \texttt{ Minimize~~} D & \label{obj_dijk} \\
\texttt{  subject to }  \nonumber \\
    n_{ft} \leq d_{ft} & \leq D, \ \forall f \ \texttt{and} \ t = a_f, \label{dijk1} \\
    d_{ft} + n_{gt} & \leq d_{gt} \ \texttt{or} \ d_{gt} + n_{ft} \leq d_{ft},\ \forall f,g \ \texttt{and} \ t \in I_f \cap I_g, \label{dijk2} \\
    d_{ft} - d_{ft-1} & \leq \texttt{max}\{0,n_{ft}-n_{ft-1}\}, \ \forall f \ \texttt{and} \ t \in(a_f,b_f], \label{dijk3} \\
    d_{ft-1} - d_{ft} & \leq \texttt{max}\{0,n_{ft-1}-n_{ft}\}, \ \forall f \ \texttt{and} \ t \in (a_f,b_f], \label{dijk4}
\end{align}

where, $D$ is the total number of counters required, $I_f$ is the check-in interval of
flight $f$ (or counter operating time), $n_{ft}$ is the number of desks required for the
check-in process of flight $f$ in period $t$ ($t \in I_f$), and $d_{ft}$ is the largest
desk number assigned to flight $f$ in period $t$ ($t \in I_f$). Constraint (\ref{dijk1})
ensures that counters assigned to flights do not exceed $D$, constraint (\ref{dijk2})
ensures that two flights are not assigned the same counter in a TW and also avoids
overlap, constraints (\ref{dijk3}) and (\ref{dijk4}) ensure that required counters are
opened or closed at the beginning of a TW. Model inputs are the counter requirements
($n_{ft}$) for flights determined using simulation. Queuing theory and simulation methods
are both evaluated, it is observed that queueing theory cannot be applied since the
arrival distribution of passengers for different flights is not homogeneous and cannot
reach a steady state. Simulation and the model ((\ref{obj_dijk})-(\ref{dijk4})) were
applied to data from a Dutch airport and it was observed that increased problem size
(increase in number of flights and planning horizon) resulted in an increase in
computation time.
 A similar mathematical model is presented by \cite{Duin2006}. The model is
adapted from RPSP (see \cite{pinedo1998operations}).

\begin{figure}
  \centering
  \includegraphics[scale=0.5]{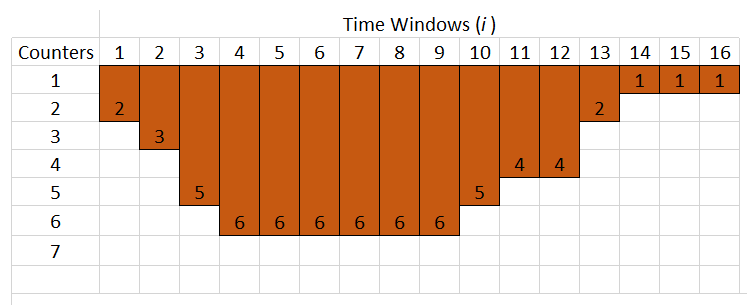}
  \caption{Task Structure for a given Counter Allocation} \label{TaskStructure}
\end{figure}

\cite{opt2019} also present an ILP for allocating adjacent counters. The ILP
((\ref{stageTwoobj})-(\ref{stageTwoconFour})) has been shown to solve problems with a
large number of departures in less time compared with formulations by
\cite{araujo2015optimizing} and \cite{bruno2010mathematical}. Tasks are defined as
time-varying counter requirement for one or more departures (see
Fig.\ref{TaskStructure}). This definition is similar to that of counter profiles by
\cite{chun1999intelligent}. The difference being that there cannot be multiple task
structures whereas the shape of a counter profile can be changed using different
operators. A counter-TW pair, $(k,\ i)$ is fixed and variables $w_{ki}$ and $ m_t$ are
defined by

\begin{eqnarray}
  w_{ki} &  =  & \sum_{t \in \mathcal{T}_i} \sum_{h=1}^{c_i(t) \wedge k} y_{t(k-h +1)}, \label{eq13:wkidef} \\
   m_t &  = & \max \{c_i(t): i_o(t)\leq i \leq i_c(t)\},  \label{eq14:mtdef}
\end{eqnarray}

where $i_o(t)$ is the counter opening time and $i_c(t)$ is the counter closing time for
task $t$, $\mathcal{T}_i$ is the set of all tasks $t$ such that $i_o(t)\leq i \leq
i_c(t)$, $w_{ki}$ is the number of tasks that use counter-TW combination $(k,\ i)$,
$y_{tk}$ is a binary variable, equal to 1 if task t starts at counter k, and $\sum_k
ky_{tk}+m_t-1$ is the largest counter number used by task $t$.

\begin{align}
 \texttt{ Minimize~~~ } s &   \label{stageTwoobj} \\
 \texttt{ subject to }& \nonumber \\
             & \sum_{t \in \mathcal{T}_i} \sum_{h=1}^{c_i(t)\wedge k} y_{t(k-h +1)} \leq 1\texttt{ for all } (i,k), \label{stageTwoconOne} \\
             & \sum_k y_{tk} = 1   \texttt{ for all } t, \label{stageTwoconTwo} \\
             & \sum_k ky_{tk} + m_t - 1 \leq s, \texttt{ for all } t, \ \label{stageTwoconThree} \\
             & y_{tk} \in \{0,\ 1\} \texttt{ for all } t,\ k, \label{stageTwoconFour}
\end{align}

Constraint (\ref{stageTwoconOne}) ensures that no counter is allocated to more than one
task in any TW. Constraint (\ref{stageTwoconTwo}) ensures that each flight is allocated
counters and constraint (\ref{stageTwoconThree}) limits the total number of counters used
to $s$.

\subsection{Some Real-life Airport Applications}
This section discusses models proposed for solving real-life airport problems. We focus
on the airport problem considered by the authors, inputs required for the models
considered, and prominent disadvantages and advantages to using these models.
\cite{joustra2001simulation} present a simulation model and \cite{dijk2006} present ILPs
in addition to simulation for check-in counter allocation, \cite{atkins2003right} present
a simulation model with input data on pre-board screening, shift scheduling, passenger
arrivals, and level of service for determining staff scheduling and resource
requirements. Simulations were run till the staff schedule obtained satisfied the level
of service at the airport. Operations in airports vary in type of check-in, type of queue
permitted, flexibility of counter usage (for business or economy class), passenger
behaviour, common or dedicated counters, check-in periods, baggage collection centers and
sorting of baggage, queuing area etc. Consequently, simulation models consider different
features and requirements of the airport under consideration.
\cite{bruno2010mathematical} present two models for counter determination. The model by
\cite{lous2011modelling} considers the baggage belt direction, amount of baggage allowed
per baggage area, the maximum number of people allowed to queue at a counter, adjacency
of counters allocated, counter location preferences by an airline, but aiming to create a
flexible model results in a complicated model, and involves a large number of additional
computations, especially in allocating preferred counters to airlines. In the study by
\cite{al2016optimization}, some counters are always kept unused in each zone of the
airport to cushion the airport operator against sudden increase in traffic. The author
overlooks the model objective function which is constant, the model is thus defective and
consequently, the objective of minimizing the counters allocated may not be achieved.
\cite{felix2017hybrid} developed a hybrid discrete-event and agent based simulation model
to assess the performance of check-in process at the airport of Lisbon. Passenger
behaviour, the sequence of tasks performed in the check-in area and the physical layout
of the check-in area are incorporated in the model. The variations in the check-in
process are attributed to variation in these aspects. The model by \cite{felix2017hybrid}
is comparable to the model by \cite{lous2011modelling} in that it considers almost the
same set of factors effecting passenger check-in.

In \cite{felix2017hybrid}, the simulation model proposed aims to explore various
scenarios at check-in area and enables airport operators to choose the best position and
assignment of different types of check-in facilities (counters, kiosks, etc) to airlines
(counters and kiosks differ in the services offered). However, in
\cite{lous2011modelling}, baggage collection centers with limited capacity and limited
queueing area are incorporated in the ILP proposed. \cite{lamphai2016} also models the
problem for an airport in Thailand. The excel based software SimQuick is used to build
the model for simulation. The arrival distribution of passengers, check-in service time
(or its distribution), type of queue, are input to the software. Different structures at
the airport such as the check-in counters, entrances, exits, queue capacity, and other
processes can be defined and average time for passenger flow through these processes
including the check-in counters can be observed. The paper concludes with different
results on the efficiency of counter allocation in terms of passenger waiting times and
queue length.

\section{Different Approaches to Counter Allocation}

This section discusses different approaches for check-in counter assignment. The
approaches used to model the problem are classified below based on the problem solved and
procedures used.

\subsection{Simulation for Counter Allocation}
Initial attempts to solve the counter allocation problem were made by simulation of
resource requirement at airports. Constraint satisfaction algorithms were presented by
\cite{chun1996scheduling} and \cite{chun1999intelligent}. Subsequently, simulation of
passenger flow at the airport terminal was proposed by \cite{wong1998development} and
passenger flow from terminal entrance to boarding was simulated by
\cite{kiran2000simulation}. \cite{wong1998development} focus on passenger traffic
characteristics and their impact on terminal operations. Simulation of resources may also
be done to identify delays at the check-in system and create scenarios that will improve
the efficiency (\cite{appelt2007simulation}). Real-time system data was used by
\cite{ros2017improving} to simulate system behaviour and test different conditions and
scenarios with the objective of optimizing the check-in procedure at Brisbane Airport.
The author presents a dynamic check-in procedure where rapid changes in the check-in
procedure are possible in real-time. Simulation was used for counter determination by
\cite{dijk2006}. For analysis of counter allocation and waiting time of passengers,
\cite{bevilacqua2010analysis} compared both waiting times and other parameters estimated
by queuing theory and simulation. \cite{chun1996scheduling} presents algorithms by
modelling the problem as a multidimensional placement problem. \cite{chun1996scheduling}
proposes a two-dimensional approximation for counter scheduling where space (counters)
and time are assumed to be part of a Gantt chart. The constraint satisfaction algorithms
presented consider all possible counter profiles (\ref{ChunCounterGlob}) for each flight.
The final version of the algorithm aims to further improve the solution by considering
different shapes per counter profile and by imposing constraints on the minimum number of
counters allocated to a flight based on the queue length restrictions. The main drawback
of this method is that final allocation depends on the order of flights chosen by the
algorithm for allocation. For the final algorithm, counter profile globs are defined as
multidimensional objects which include additional dimensions such as queue lengths,
waiting time and baggage restrictions. Also, in case all the flights cannot be allocated
to counters, despite re-shaping the counter profile globs and minimizing the counters
allocated, the algorithm backtracks, de-assigns the most recent allocation and proceeds
to allocate for another flight. This is a major drawback for airports with a large number
of departures. There may be too many backtracks in deriving the final solution. Also, the
final allocation depends on the order of the flights chosen. Choosing the best counter
profile has not been addressed in this paper,
which is very much needed for arriving at a good decision.  \\
\cite{chun1999intelligent} present an intelligent resource simulation system (IRSS). IRSS
predicts the check-in counter requirement at an airport. This is a software which takes
 the airport flight data, passenger turn up data, check-in
counter data, and other model parameters as input to generate a simulation model. The
IRSS proposed also has a graphic user interface to simulate and animate the check-in
counter queueing for a single flight. Simulation parameters such as tolerable passenger
waiting time, queue lengths, service rates, check-in times for flights, the number of
check-in counters, number of passengers or time of departure can be changed to match the
reality at the airport and to observe the change in counter allocation. The objective is
to find a counter profile (a counter assignment solution) which saves the maximum counter
time compared to the current counter allocation and maintains the desired service level
quality. The authors also account for the stochastic processes such as arrival rates. The
Check-in Counter Allocation System (CCAS) presented by \cite{chun1996scheduling} together
with the IRSS are capable of addressing check-in counter allocation problem at an airport
and were used together in the Kai Tak International Airport. The solutions obtained are
evaluated by analysing the waiting times of passengers and queue lengths. Though the
algorithm tries to achieve the best possible solution, a major disadvantage is that not
all possible counter profiles can be simulated. The authors examine some %(five to eight)
counter profiles to arrive at the best counter allocation (\cite{chun1996scheduling}).
Due to this we cannot ensure that the service quality level provided is the best
possible. Also, for a large number of flights and time windows, algorithms with
simulation will take a large amount of time.  \\

\cite{bevilacqua2010analysis} present a case study and an analysis of counter allocation
through simulation. Their objective is to calculate minimum counter requirement of each
flight such that average waiting time does not exceed the maximum allowable limit.
Passenger arrivals and counter operations are simulated for a single flight or a group of
flights. Poisson distribution was used to generate arrivals. The model assumes steady
state of the system. In the case study presented, counter operations are modelled as a
queueing theory system and stability conditions for common check-in were calculated.
Simulation analysis shows that common check-in is better than dedicated check-in.
Simulation results are compared with standard results from queuing theory. The effect of
varying arrival time distribution and number of counters allocated to a flight on the
average waiting time is computed. It is observed that queueing theory is more applicable
for common check-in and simulation is more suited for dedicated check-in.

Using simulation for decision making process has the main disadvantage that it explores
only a small subset of the whole possible scenarios that can be reached by the system
under study, thus reducing its optimization potential (\cite{mota2015check}). To overcome
this problem many authors have proposed mathematical models for determining counter
requirements and for ensuring adjacency in counter allocation. Some of the models are
discussed in the sections ~\ref{SectionCDetermination} and ~\ref{SectionCAllocation}.

\subsection{Network Model for Counter Allocation}
\begin{figure}
\centering
\begin{subfigure}{.49\textwidth}
  \centering
  \includegraphics[scale=.8]{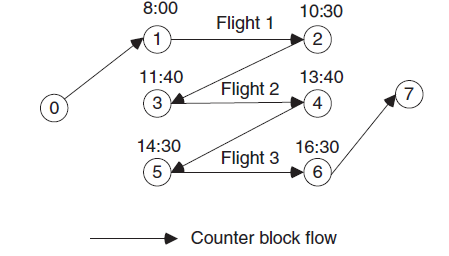}
  \caption{Counter Block Flow}
  \label{fig:counterblockflow}
\end{subfigure}
\begin{subfigure}{.49\textwidth}
  \centering
  \includegraphics[scale=.6]{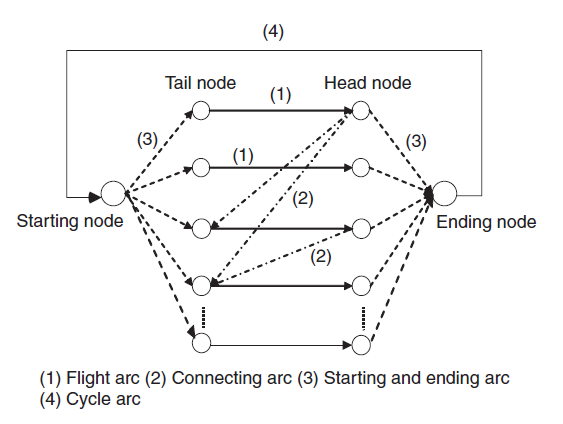}
  \caption{Multiple Counter Block Flow Networks}
  \label{fig:multiplecounterblockflows}
\end{subfigure}
\caption{Networks for Counter Allocation by \cite{tang2010network}}
\label{network}
\end{figure}

\cite{tang2010network} developed a network flow model for allocating blocks of counters
to airlines. The network flow model aims to optimize counter usage at a Taiwan airport.
The model uses predefined blocks of counters at the airport and opening and closing
counter times of flights. Each flight is defined using an arc, the opening time and the
closing time form connecting nodes. Arcs are used to denote sequential flight allocations
in a counter block (Fig.~\ref{fig:counterblockflow}).

The opening time, closing time and counter requirements for a flight are inputs to the
model, as determined by airlines. An ILP meeting all the constraints at the airport is
proposed. Flights, set of nodes and arcs that can be assigned to a block are all inputs
to the model resulting in a counter block flow network for each block. A network can also
be used to represent all the possible flight assignments to a block in a day
(Fig.~\ref{fig:multiplecounterblockflows}). The main advantages of the network flow
method are allocation of adjacent counters and convenient representation of the flight
sequence allocated to a counter block. The model involves preprocessing to eliminate
redundant constraints (in the example presented by \cite{tang2010network}, about 30\% of
the constraints are redundant) and computing parameters before constructing the model.
Due to this preprocessing, running time varies exponentially with input size (the number
of flights scheduled and the number of counter blocks allowed). The model is used to find
near-optimal solutions at the Taiwan Taoyuan International Airport, though implementation
becomes extremely complicated for assigning time-varying counters to flights. The main
drawback of their model is that the counters are predivided into blocks for allocating to
airlines which results in a lot of preprocessing. Since multiple blocks can be built with
the same counter, the number of possible blocks can be very large complicating the
process of counter allocation.

\subsection{Evolutionary Algorithms and Counter Allocation}
Some authors have proposed genetic and evolutionary algorithms for the check-in counter
allocation problem (\cite{yeung1995check}, \cite{mota2015allocation},
\cite{mota2015check}, \cite{motaandzuniga2013}). Evolutionary techniques are a group of
methods inspired by common evolutionary processes. These techniques are expected to
provide good solutions, i.e. solutions that are close to optimal but may not be
optimal(see \cite{Goldberg:1989:GAS:534133} and \cite{mota2015check}). The efficiency of
these techniques relies highly on parameters that drive the selection procedure (see
\cite{mota2015check} and \cite{affenzeller2009genetic}). Genetic algorithms are a part of
evolutionary techniques. Genetic algorithms (GAs) are defined as efficient, adaptive and
robust search and optimization processes that are applied in large and complex search
spaces. GAs are modelled on the principles of natural genetic systems where the genetic
information of each individual or potential solution is encoded in structures called
chromosomes. GAs compute a fitness function for directing search in more promising areas.
Each individual has an associated fitness value, which indicates its degree of goodness
with respect to the solution it represents. GAs search from a set of points called a
population and various biologically inspired operators like selection, crossover and
mutation are applied to obtain better solutions (\cite{bandyopadhyay2007classification}).

\cite{yeung1995check} use fitness directed scheduling to develop an airport check-in
counter allocation system based on genetic algorithms. Populations of individuals are
genetically bred according to Darwinian principles, i.e., reproduction of the fittest and
crossover operations. Each individual represents a check-in counter allocation plan for
one day (it is the allocation on a Gantt Chart). Each individual also has an associated
fitness measure. Fitness measure is in terms of the number of overlaps found in an
allocation plan. Lesser the overlaps, fitter the allocation plan. The fittest individual
is the best allocation plan for that day. A population of individuals is randomly created
and fittest individuals are selected for the crossover operation. The crossover operation
is to create offspring counter allocation plans from the selected individuals in the
population. By recombining randomly chosen assignments of the fittest allocation plans,
we produce new allocation plans. This new population replaces the previous population and
the entire process is repeated to create new generations. The best allocation plan that
appeared in any generation is the best plan for check-in counter allocation problem.

\begin{figure}
  \centering
  \includegraphics[scale=0.5]{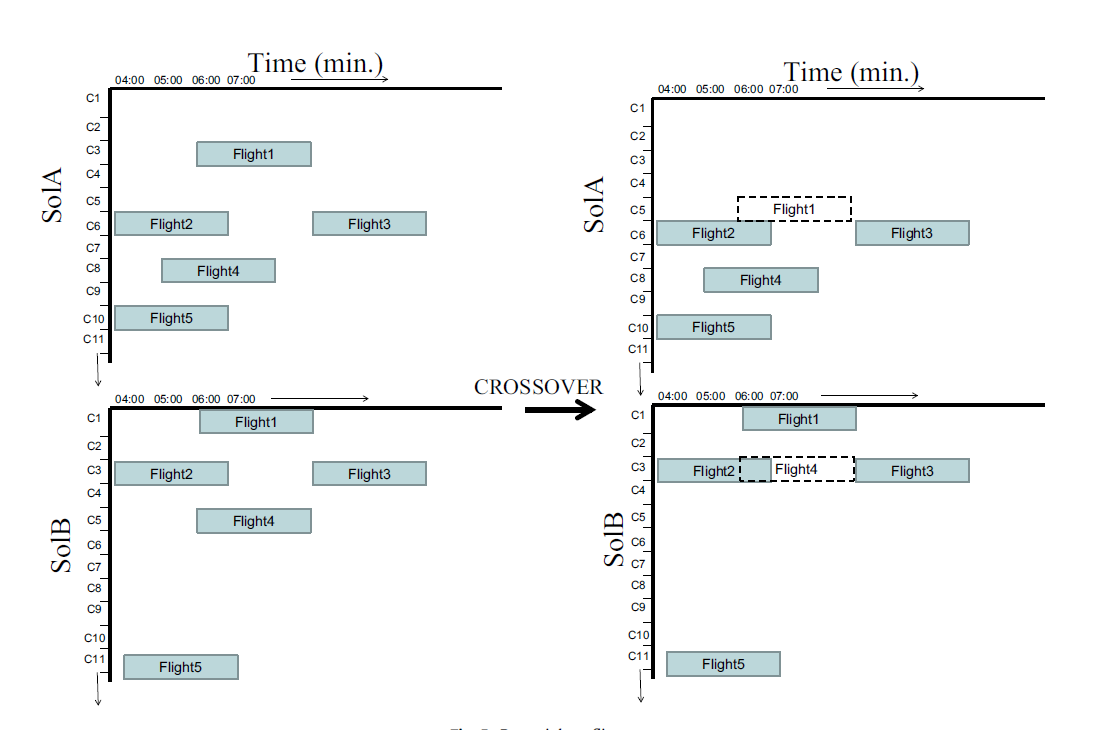}
  \caption{Crossover of two solutions} \label{crossover}
\end{figure}

\cite{mota2015check} presents a methodology that combines evolutionary techniques and
simulation and aims to provide a solution which is better than the solution obtained by
applying these techniques independently. The algorithm uses a brute-force approach.
Flights are allocated sequentially, taking into account all the constraints in flight
allocation such as no overlap, counters opened three hours before departure etc. After
all the flights are allocated an initial solution is obtained. This is similar to the
algorithm first-fit. In order to get varying solutions, the flight order is changed
before each allocation. A population of Counter allocation plans is thus obtained. Next,
the solutions are converted into vectors (chromosomes) with information that will be used
by the evolutionary algorithm. Crossover operations are performed (see
Fig.\ref{crossover}) to improve the existing solutions such that feasibility of the
generated solution is maintained. To perform a crossover between any two initial
solutions, flights are selected randomly from each of the two solutions and compared. For
a pre-decided percentage of flights, with matching counter opening and closing times,
counter locations are interchanged. The crossover procedure is performed on solutions
with higher fitness values. The resulting solutions with a high measure of fitness are
then retained and again crossed over. Though feasibility of the resulting solution is
ensured, the computation time to determine feasibility for crossover is very large and
increases with the size of the input.  For instance, for an airport with about 2500
departures in a week, to perform a crossover, 2500*2500 flight combinations need to be
checked for possibility of a crossover. The cost function is computed for each generation
and checked for improvement. Since the problem is multi-objective, an objective function
is computed as a fitness measure. In the algorithm proposed, the solutions are improved
till a stop condition is reached. The stop condition is arbitrarily determined and
feasible solutions obtained this way are analysed in real-life conditions using
simulation.

A major problem with genetic and evolutionary techniques is that the solution(s) with
best fitness may not be anywhere near optimal, but is only relatively better among the
allocations generated. Larger the initial population of
counter assignments higher the possibility of improving the initial solution. For $n$
flights, the algorithm by \cite{yeung1995check} needs $(1 + \binom{2^n}{2}) *(2^n)$
steps, where $2^n$ is the number of possible recombinations of two counter allocations.
This results in unnecessarily prolonged and time consuming calculations for eliminating
poor solutions, thus, the number of steps for the algorithm increase exponentially
with $n$. %(Performance = O($2^n)$)

\subsection{Queuing Theory and Counter Allocation} \cite{parlar2008dynamic} propose a
dynamic programming technique based on queuing theory for determining the counter
requirement of a flight. Queueing theory results are used to model the problem. A pure
death process is used to model the arrivals. An exponential service time distribution is
assumed (Erlang distribution is also explored). The rate of service is observed as
proportional to the queue size, hence, service rate is assumed to be state-dependent. The
exact forms of distributions of the arrival rate and the service rate are found. These
are used to calculate the expected number of passengers in the system and the cost of
passengers waiting to be served. Since the arrival rate is time dependent and the arrival
process is non-stationary, arrivals are observed and counter operating time is divided
into smaller subintervals with constant arrival rate. The arrival rate is then estimated.
A dynamic programming model is then used to determine counters to be opened. Counter
opening or closing decisions are expected to be made every 20 minutes to minimize the
total expected cost. \cite{parlar2008dynamic} propose a model to determine optimal time
varying counter allocation for a single flight. \cite{parlar2013allocation} propose
static counter allocation policy for a single flight. Their objective is to minimize the
expected total cost of waiting, counter operation, and passenger delay which the authors
show to be convex in the number of counters allocated. The authors also introduce a
service level constraint to ensure that a certain percentage of passengers are served in
each TW.

\section{Related Scheduling Problems}

Two problems related to ACAP, the two dimensional strip packing problem and resource
constrained project scheduling (RPSP), have been studied extensively (see
\cite{lodi2002two}, \cite{amoura2002scheduling}, \cite{blazewicz1986scheduling},
\cite{Duin2006} for further details). The two dimensional strip packing problem comprises
of allocating rectangular items to a larger standard size rectangle with the objective of
minimizing waste. The problem, $P|fix_j|C_{max}$ in multiprocessor scheduling, after
swapping time and place is equivalent to the adjacent resource allocation problem with
rectangular units (for further details see \cite{Duin2006} and
\cite{amoura1997scheduling}). The problem with irregularly shaped units, such as
polyominoes (see Fig.\ref{VanDijk1}), is similar to RPSP. Hence, a well known integer
programming formulation of RPSP (see \cite{pinedo1998operations}, \cite{Duin2006}) can be
modified to solve ACAP (\cite{Duin2006}).

\cite{Duin2006} discusses similarities of the counter allocation problem with other
resource allocation problems in literature such as the two-dimensional strip packing
problem (see \cite{lodi2002two}) and the resource constrained project scheduling problem
(see \cite{blazewicz1986scheduling} and \cite{du1989complexity}).

Staff rostering/human resource management problems at the airport check-in counters are
discussed by \cite{lin2015ground}, \cite{zamorano2018task}, \cite{rodivc2017airport} and
\cite{xin2014design}. \cite{xin2014design} discusses both counter allocation and staff
rostering.

\cite{bruno2018decision} propose a model to optimize shift scheduling decisions of desk
operators and service level measured in terms of passenger waiting times at the counters.
A real-life case study has been presented for two airports in Italy.

\cite{hsu2005scheduling} study optimal facility purchase and replacement.
\cite{brunetta1999operations} evaluates an airport terminal and estimate delays due to
facilities such as the check-in counters. \cite{fayez2008managing} estimates the
passenger flow through the airport. Efficient use of airport capacity is discussed by
\cite{pacheco2003managerial}.

\cite{yan2008reassignments} and \cite{yan2014common} have worked on reassignments of
counter allocations in case of sudden unexpected events at the airport such as change in
flight schedule, baggage belt malfunction, airport closure and other disturbances to the
planned counter allocation. A mathematical formulation has been presented by
\cite{yan2014common} with the objective of reducing the impact of unforeseen
circumstances at the airport. An inconsistency for a flight is defined as the deviation
between original and reassigned counters. The model aims to reduce the inconsistencies in
assignment. A heuristic is proposed for solving the model. The model is modified into two
relaxations to obtain an upper bound and lower bound. The two relaxations of the model
are repeatedly solved till the difference between the two is lower than a predefined
limit. The main advantage of the formulation is that it helps the airport authorities
with reassignment of counters to restore normalcy and contain the impact of a disturbance
to as few flights as possible. Some numerical validation to the model is given. A major
disadvantage is that all flights may not be reassigned adjacent counters. Different
counters in different TWs may result in confusion for the customers, as it is not
possible to shift passengers in a queue from one counter to another.

\section{Conclusion}

Various methods of solving the counter allocation problem at airports are presented in
this paper. It is observed that determining counter requirements for flights and then
allocating adjacent counter space is most suitable to obtain a practical solution to the
problem.

\bibliography{CAPbibfile}
\bibliographystyle{ecta}
%\printbibliography %[title=Bibliography]
\end{document}